# Value Function Gradient Learning for Large-Scale Multistage Stochastic Programming Problems

Jinkyu Lee[a], Sanghyeon Bae[a], Woo Chang Kim[a,*], Yongjae Lee[b,*]

**Abstract**

A stagewise decomposition algorithm called "value function gradient learning" (VFGL) is proposed for large-scale multistage stochastic convex programs. VFGL finds the parameter values that best fit the gradient of the value function within a given parametric family. Widely used decomposition algorithms for multistage stochastic programming, such as stochastic dual dynamic programming (SDDP), approximate the value function by adding linear subgradient cuts at each iteration. Although this approach has been successful for linear problems, nonlinear problems may suffer from the increasing size of each subproblem as the iteration proceeds. On the other hand, VFGL has a fixed number of parameters; thus, the size of the subproblems remains constant throughout the iteration. Furthermore, VFGL can learn the parameters by means of stochastic gradient descent, which means that it can be easily parallelized and does not require a scenario tree approximation of the underlying uncertainties. VFGL was compared with a deterministic equivalent formulation of the multistage stochastic programming problem and SDDP approaches for three illustrative examples: production planning, hydrothermal generation, and the lifetime financial planning problem. Numerical examples show that VFGL generates high-quality solutions and is computationally efficient.

[a] Department of Industrial and Systems Engineering, Korea Advanced Institute of Science and Technology (KAIST)

[b] Department of Industrial Engineering, Ulsan National Institute of Science and Technology (UNIST)

[*] Corresponding authors: wkim@kaist.ac.kr, yongjaelee@unist.ac.kr

# 1. Introduction

Multistage stochastic programming (MSP) problems arise in a broad range of areas where decisions should be made under uncertain environments. Such optimization problems include capacity expansion planning, asset liability management, and hydropower production planning (Shiina & Birge, 2003; Cariño et al., 1994; Fleten & Kristoffersen, 2008). The stochastic process, which models the uncertain environment, is generally approximated by a set of finite scenarios, which is called a "scenario tree." Then, the problem can be replaced by a single large deterministic equivalent problem. However, the deterministic equivalent problem often becomes computationally intractable because the number of scenarios grows exponentially with respect to the number of stages. This "curse of dimensionality" is a critical issue for MSP problems with long-term planning horizons.

In this study, a novel stagewise decomposition algorithm called "value function gradient learning" (VFGL) is proposed. It approximates value functions as parametric convex functions and learns parameters that closely approximate the gradient of the true value functions. To be more specific, the proposed algorithm minimizes the stagewise estimation error of the gradient of the approximated value function based on stochastic gradient descent to find the best parameters within a given parametric form. We show under mild regularity conditions that the approximated solution converges to the optimal solution. Furthermore, we introduce a metric called KKT deviation to measure the suitability of the given parametrization. In our numerical examples, we show that VFGL can efficiently find good quality solutions with the appropriate parameterization of the value function.

There are three unique characteristics of VFGL compared with the traditional stagewise decomposition method. First, VFGL approximates the value function using a fixed parametric functional form. The main advantage of this value function approximation approach is that the size of the subproblems remains constant for every iteration. However, the piecewise linear approximation of the value function results in an increasing size of the subproblems in every iteration because of the increasing number of linear cuts. This makes VFGL computationally suitable for stochastic programs that require many iterations for the value function approximation.

Second, VFGL directly samples a random process from its distribution without the scenario tree approximation. A scenario tree approximation of the distribution of continuous random variables is a nontrivial task. There are many reports in which various scenario tree generation algorithms are discussed (Pflug, 2001; Gülpınar et. al., 2004; Heitsch & Römisch, 2009). However, they are rather heuristic, and the choice of the scenario tree generation algorithm for a particular problem is an open-ended question. Therefore, MSP problems are often solved multiple times under various scenario trees to check the distributional approximation error indirectly. VFGL, by contrast, is free of this issue. However, VFGL might involve the extra step of verifying the suitability of the chosen parametric form of the value function for a specific problem.

Finally, VFGL is an online learning algorithm that learns the parameters online during the exploration of feasible paths. A parallel computation is easy for VFGL because forward explorations can be computed in parallel. An iteration of SDDP involves a forward pass that samples stagewise trial solutions and backward pass that solves all possible stagewise subproblems in a reverse stagewise order. However, an iteration of VFGL involves only forward simulation to sample stagewise trial solutions and update the value function by adjusting parameters in a stagewise order. Therefore, each iteration of VFGL tends to be much faster than that of the traditional stagewise decomposition algorithms. However, one drawback of online learning is that the value function update information is passed on to only the prior stage. Therefore, the number of iterations should be greater than the number of stages so that the final stage information is passed onto the first-stage value function

Three illustrative examples are provided to compare the MSP, SDDP, and VFGL approaches: production optimization, a hydrothermal scheduling problem, and a discretized lifetime portfolio optimization problem. The numerical experiment clearly shows the numerical potential of VFGL compared with SDDP, which is considered a state-of-the-art stagewise decomposition algorithm. The numerical study further shows that the parametric value function form can be recycled to a limited extent for the perturbation problem depending on the problem and perturbation. The performance evaluation of VFGL includes a comparison of various parametric value function forms, where the usefulness of KKT deviation is illustrated.

The remainder of this article is organized as follows. Section 2 reviews literatures on various methods for sequential decision making problems under uncertainty to clarify the difference between the VFGL and existing algorithms. In Section 3, the considered problem is mathematically defined, and the VFGL is derived. Three illustrative examples are shown in Section 4. The conclusions are provided in Section 5.

## 2. Literature review

Stochastic optimization encompasses a broad range of problems involving optimal decision making over a period in an uncertain environment. Powell (2019) classified diverse studies of stochastic optimization into fifteen categories in terms of key modelling characteristics, including the problem statement regarding optimality, state variable, interaction between uncertainty and decision variables, system dynamics modelling, and objective functions. In this study, we focused on discrete time sequential stochastic optimization models, which discretize time to achieve finite or countable time transitions. Powell (2019) classified discrete time models into the following two major streams. 1) From the operations research community. Herein, the multistage stochastic programming (MSP) is developed from deterministic optimization; subsequently, the MSP is extended to Markov decision process (MDP) to handle large-scale problems. 2) From the computer science community. Herein, MDP models are directly addressed using reinforcement learning (RL).

2.1. Multi-stage stochastic programming (MSP) approach

The MSP framework deals with the optimization of sequential linear/convex/nonconvex programming problems, where parameters of each program depend on the realization of random variables and decisions from the prior stage(s). An MSP problem is solved approximately using its deterministic equivalent problem, which is constructed by a finite scenario tree that represents the evolution of the underlying stochastic process. However, it often becomes computationally intractable because the number of scenarios increases exponentially with respect to the number of stages and/or the number of nodes per stage. To manage this "curse of dimensionality", various decomposition-based algorithms have been proposed for solving large-scale stochastic programs.

2.1.1. Scenario decomposition methods

Scenario decomposition methods decompose stochastic programs in a scenario-wise manner by relaxing the nonanticipativity constraints, which prevent the use of information unavailable at the time. The two most popular examples for this are the dual decomposition algorithm (Carøe & Schultz, 1999) and progressive hedging algorithm (Rockafellar & Wets, 1991).

The dual decomposition algorithm was proposed to solve linear multistage stochastic integer programming problems. It relaxes the nonanticipativity constraints to derive the Lagrangian dual problem and employs a branch-and-bound method to obtain the optimal solution using the lower bound information obtained from the dual problem. In contrast, the progressive hedging algorithm introduces modified scenario-wise subproblems by relaxing the nonanticipativity constraints through an augmented Lagrangian function that incorporates a multiplier term and quadratic penalty terms. The algorithm iteratively adjusts the solutions to determine implementable and admissible solutions that satisfy the nonanticipativity constraints and lie within the scenario-wise feasible region.[1]

However, scenario decomposition algorithms visit all scenarios during every iteration; otherwise, the solutions may violate the nonanticipativity constraints. Therefore, the scenario decomposition approach may be inappropriate for MSP problems with extremely large scenario trees.

2.1.2. Stagewise decomposition methods

Stagewise decomposition algorithms break down the stochastic program into stagewise subproblems and subsequently determine solutions by sequentially solving the subproblems. The objective function of each subproblem incorporates a special function called a "value function" (often referred to as a recourse function or a cost-to-go function), which reflects the future consequences of an immediate decision.

The stochastic dual dynamic programming (SDDP) algorithm, which is a sampling-based variant of a nested Benders decomposition, was introduced by Pereira and Pinto (1991). SDDP iterates between the forward and backward passes. In a forward pass, the algorithm samples a bundle of scenarios to obtain feasible policies, whereas in a backward pass, the value functions are improved by constructing subgradient cuts around each of the stagewise solutions of the feasible policies.[2] However, the size of the stagewise subproblems increases as the algorithm proceeds because the number of subgradient cuts increases with each iteration. This is one of the major drawbacks of SDDP.

---

[1] Rockafellar and Wets (1991) showed that the progressive hedging algorithm converges to the optimal value for convex optimization problems. For nonconvex cases, although convergence is not guaranteed, the algorithm is shown to produce high-quality solutions (Watson et al., 2011).

The stagewise decomposition approach is analogous to addressing MSP problems under the MDP context. A stagewise decomposition method finds an optimal policy under Markovian state transitions, which significantly reduces computational costs. However, state transitions in scenario trees of MSP are often not Markovian because they can depend on the historical state path. Therefore, the stagewise decomposition approach incorporates Markovian state transitions by employing scenario lattices instead of scenario trees. Essentially, for any node, its transition probability distribution is shared among all nodes in the same stage; thus, a value function can be shared among nodes within the stage.

2.2. Reinforcement learning (RL) approach

Recently, reinforcement learning (RL) algorithms have gained considerable attention as a solution to MDP problems. RL algorithms use a data-driven approach; essentially, they find an optimal policy with trial-and-error under the simulated environment. In particular, numerous RL algorithms fix the parametric form of value/policy functions because it is advantageous for learning from large, simulated data sets (Geist and Pietquin, 2013; Mnih et al., 2013; Mnih et al., 2016). In addition to the development of deep learning models, some RL models have exhibited impressive performances that overwhelm human experts in numerous applications, such as Go and video games. However, they face difficulties in handling continuous action spaces with complex state-dependent restrictions.

Although both RL models and stagewise decomposition algorithms (for MSP) solve sequential stochastic optimization problems under the MDP framework, they exhibit the following two distinct differences. First, MSP considers continuous action spaces that are constrained by linear equalities and convex inequalities. Furthermore, these constraints are dependent on the current state. However, RL usually considers finite or locally unconstrained action spaces. Second, immediate reward structure and state transition dynamics are explicitly known in MSP. Hence, MSP problems can be directly solvable using optimization solvers. On the other hand, RL models employ the trial-and-error approach, which would require considerable amount of data and time.

3. **Value function gradient learning**

---

[2] The convergence of SDDP was shown for risk-neutral multi-stage stochastic linear programs by Chen et al. (1999) and Philpott et al. (2008) and for risk-averse multi-stage stochastic linear programs by Guigues et al. (2012). Recent studies applied SDDP to multistage stochastic convex programs and investigated its convergence for risk-neutral and risk-averse cases (Girardeau et al., 2015; Guigues, 2016).

The purpose of the value function gradient learning (VFGL) algorithm is to find optimal policies for large-scale MSP problems. We tackle the problem with a stagewise decomposition approach, with the main focus on approximating the value functions efficiently. The SDDP approximates value functions using piece-wise linear functions, whereas VFGL approximates value functions with parametric functions as in numerous RL algorithms. However, the parameter learning process of VFGL is very different from those of RL algorithms. RL algorithms sample costs/rewards based on trial-and-error, and parameters are updated from those cost/reward signals. In contrast, VFGL samples actions based on simulation, but it learns the parameters using the stagewise gradient information of the true value function by fully exploiting the duality theory in optimization.

In this section, we present the value function gradient learning (VFGL) algorithm, which is a new stagewise decomposition algorithm for large-scale multistage stochastic programming problems. First, we define the problem setting for stagewise decomposition algorithms (including VFGL and SDDP) in Section 3.1. In Section 3.2, we describe the main idea of VFGL and why it will work. Section 3.3 details the learning procedure of the optimal parameters based on a stochastic gradient descent. Finally, we discuss how to choose the appropriate parameterization in Section 3.4.

## 3.1. Problem setting

A sequence of decision-making over multiple periods is considered, where the relevant stochastic process is gradually realized. There is a stochastic process $\xi_{[T]} = (\xi_1, \ldots, \xi_T)$, with $\xi_1$ deterministic, and a decision process $x_{[T]} = (x_1, \ldots, x_T)$, where $x_t \in \mathbb{R}^{n_t}$ in which $n_t$ is the number of dimensions of the decision variable at stage $t$. It is further assumed that each $\xi_t$ has a finite moment. Each decision was made using only the information available at that time. In particular, the following sequence of decisions and observations is assumed (Shapiro et al., 2009):

$$decision(x_1) \rightsquigarrow observation(\xi_2) \rightsquigarrow decision(x_2) \ldots \rightsquigarrow observation(\xi_T) \rightsquigarrow decision(x_T)$$

Here, the decision process $x_t$ is $\mathcal{F}_t$-measurable, where $\mathcal{F}_t$ is the sigma algebra generated by a stochastic process $(\xi_1, \ldots, \xi_t)$. A $T$-stage stochastic program in nested form is formulated as follows:

$$\min_{x_1 \in \mathcal{X}_1} f_1(x_1) + \mathbb{E}\left[\min_{x_2 \in \mathcal{X}_2(x_1,\xi_2)} f_2(x_2, \xi_2) + \mathbb{E}_{\cdot|\xi_{[2]}}\left[\ldots + \mathbb{E}_{\cdot|\xi_{[T-1]}}\left[\min_{x_T \in \mathcal{X}_T(x_{T-1},\xi_T)} f_T(x_T, \xi_T)\right]\right]\right] \quad (1)$$

where $\mathbb{E}_{\cdot|\xi_{[t]}}$ is a conditional expectation operator with respect to $\xi_{[t]}$, the history of the data process up to stage $t$. For $t = 1$, $f_1(x_1)$ is a deterministic objective function that is convex in $x_1$, and $\mathcal{X}_1 = \{x_1 : g_{1,i}(x_1) \leq -h_{1,i}, i = 1, \ldots, p_1 \text{ and } l_{1,j}(x_1) = b_{1,j}, j = 1, \ldots, q_1\}$ is defined by a deterministic convex feasible region for $x_1$, where $g_{1,i}$ is a twice-differentiable convex function in $x_1$, $l_{1,j}$ is a linear function in $x_1$, and $h_{1,i}$ and $b_{1,j}$ are deterministic coefficients. For $t = 2, \ldots, T$, $f_t(x_t, \xi_t)$ is an objective function that is convex in $x_t$ and dependent on $\xi_t$, and $\chi_t(x_{t-1}, \xi_t)$ is a convex feasible region for $x_t$ given $x_{t-1}$ and $\xi_t$. In particular, $\chi_t(x_{t-1}, \xi_t) := \{x_t : g_{t,i}(x_t, \xi_t) \leq -h_{t,i}(x_{t-1}, \xi_t), i = 1, \ldots p_t$ and $l_{t,j}(x_t, \xi_t) = b_{t,j}(x_{t-1}, \xi_t), j = 1, \ldots, q_t\}$, where $g_{t,i}$ is a twice-differentiable convex function in $x_t$, $l_{t,j}$ is a linear function in $x_t$, $h_{t,i}$ is a twice-differentiable convex function in $x_{t-1}$, and $b_{t,j}$ is a linear function in $x_{t-1}$.

The usual MSP approach solves Eq. (1) in two steps: 1) approximating the distribution of the stochastic process $\xi_{[T]}$ by a scenario tree with a finite number of realizations, and 2) solving a large deterministic equivalent convex optimization problem under the approximated scenario tree. However, the deterministic equivalent problem often becomes intractable; hence, the problem is decomposed into subproblems in a stagewise manner.

In this study, the following assumptions are made, which are common in stagewise decomposition approaches, such as SDDP, for sequential decision-making problems.

*(A1) Stagewise independence*: For each stage $t = 1, \ldots, T$, $\xi_t$ is independent of $\xi_{[t-1]}$.

*(A2) Relatively Complete Recourse*: For each stage $t = 1, \ldots, T$, for any feasible $x_{t-1}$, and for any realization $\xi_t^s$ of $\xi_t$, $\chi_t(x_{t-1}, \xi_t^s)$ is bounded and nonempty.

Under (A1), the nested structure of Eq. (1) yields the following Bellman equation.

For $t = T, \ldots, 2$

$$\begin{aligned} \mathcal{V}_t(x_{t-1}, \xi_t) &= \inf_{x_t \in \mathcal{X}_t(x_{t-1}, \xi_t)} \{f_t(x_t, \xi_t) + V_{t+1}(x_t)\} \\ V_t(x_{t-1}) &:= \mathbb{E}[\mathcal{V}_t(x_{t-1}, \xi_t)] \end{aligned} \tag{2}$$

with $V_{T+1} \equiv 0$.

Under the Bellman equation Eq. (2), the optimal policy is obtained by solving the following convex problem.

*Problem (1)*

For $t = 1$,

$$x_1^* \in \arg\min_{x_1 \in \mathcal{X}_1} f_1(x_1) + V_2(x_1)$$

For $t = 2, \ldots, T$.

$$x_t^* \in \arg\min_{x_t \in \mathcal{X}_t(x_{t-1}, \xi_t)} f_t(x_t, \xi_t) + V_{t+1}(x_t)$$

with $x_t^*$ being a function of $x_{t-1}$ and $\xi_t$

## 3.2. Loss function

The main idea of VFGL is to approximate the value function $V_t$ with a fixed convex parametric function $\widehat{V}_t$ with parameter $\theta_t \in \mathbb{R}^{m_t}$. Hence, VFGL solves a sequential decision-making problem by learning the parameters $\theta_t, t = 2, \ldots, T$. The revised decomposed subproblems can be described as follows.

The current policy is obtained by solving the following convex problem.

*Problem (2)*

For $t = 1$,

$$x_1^* \in \arg\min_{x_1 \in \mathcal{X}_1} f_1(x_1) + \widehat{V}_2(x_1; \theta_2)$$

For $t = 2, \ldots, T$,

$$x_t^* \in \arg\min_{x_t \in \mathcal{X}_t(x_{t-1}, \xi_t)} f_t(x_t, \xi_t) + \widehat{V}_{t+1}(x_t; \theta_{t+1})$$

with $x_t^*$ being a function of $x_{t-1}$ and $\xi_t$

The KKT optimality conditions of the stage $t$ subproblem for $t = 2, \ldots, T$ of *Problem* (2) are given below.

*Stationarity:*

$$\nabla_{x_t} f_t(x_t, \xi_t) + \nabla_{x_t} \widehat{V}_{t+1}(x_t; \theta_{t+1}) + \sum_{i=1}^{p_t} \mu_{t,i} \nabla_{x_t}\left(g_{t,i}(x_t, \xi_t) + h_{t,i}(x_{t-1}, \xi_t)\right) + \sum_{j=1}^{q_t} \lambda_{t,i} \nabla_{x_t}\left(l_{t,j}(x_t, \xi_t) - b_{t,j}(x_{t-1}, \xi_t)\right) = 0$$

*Primal Feasibility:*

$$g_{t,i}(x_t, \xi_t) \leq -h_{t,i}(x_{t-1}, \xi_t) \qquad i = 1, \ldots, p_t$$
$$l_{t,j}(x_t, \xi_t) = b_{t,j}(x_{t-1}, \xi_t) \qquad j = 1, \ldots, q_t$$

*Dual Feasibility:*

$$\mu_{t,i} \geq 0, \; i = 1, \ldots, k$$

*Complementary Slackness:*

$$\mu_{t,i} \nabla_{x_t}\left(g_{t,i}(x_t, \xi_t) + h_{t,i}(x_{t-1}, \xi_t)\right) = 0, \; i = 1, \ldots, k_t$$

Let $(\hat{x}_t^*, \hat{\mu}_t^*, \hat{\lambda}_t^*)$ be the optimal primal-dual triple that satisfies the above KKT conditions of the stage $t$ subproblem in *Problem* (2). Then, this optimal triple of the approximated problem satisfies Proposition 1.

**Proposition 1** (*perturbed KKT condition for a subproblem with approximated value function*)

Let $(\hat{x}_t^*, \hat{\mu}_t^*, \hat{\lambda}_t^*)$ be an optimal primal-dual triple that satisfies the KKT condition of the stage $t$ subproblem of *Problem* (2). Then, $(\hat{x}_t^*, \hat{\mu}_t^*, \hat{\lambda}_t^*)$ satisfies the primal feasibility, dual feasibility, and complementary slackness of the KKT condition of the stage $t$ subproblem of *Problem* (1), while the stationarity condition is perturbed as follows.

$$\nabla_{x_t} f_t(\hat{x}_t^*, \xi_t) + \nabla_{x_t} V_{t+1}(\hat{x}_t^*) + \sum_{i=1}^{p_t} \hat{\mu}_{t,i}^* \nabla_{x_t}\left(g_{t,i}(\hat{x}_t^*, \xi_t) + h_{t,i}(x_{t-1}, \xi_t)\right)$$

$$+ \sum_{j=1}^{q_t} \hat{\lambda}_{t,j}^* \nabla_{x_t}\left(l_{t,j}(\hat{x}_t^*, \xi_t) - b_{t,j}(x_{t-1}, \xi_t)\right) = \nabla_{x_t} V_{t+1}(\hat{x}_t^*) - \nabla_{x_t} \hat{V}_{t+1}(\hat{x}_t^*; \theta_{t+1})$$

**Proof.** See Appendix A.

The KKT perturbation from the value function estimation error described in Proposition 1 is denoted with $D_t(\theta_{t+1}; x_t)$ as follows.

$$D_t(\theta_{t+1}; x_t) := \|\nabla_{x_t} V_{t+1}(x_t) - \nabla_{x_t} \hat{V}_{t+1}(x_t; \theta_{t+1})\|$$

Herein, $\|\cdot\|$ indicates the Euclidean norm, and $dist(\cdot,\cdot)$ refers to the distance measure induced by the Euclidean norm. Let $(x_t^*, \mu_t^*, \lambda_t^*)$ be the true optimal solution under the KKT optimality conditions of *Problem* (1). It is assumed that the second-order sufficiency condition (SOSC) is satisfied for $(x_t^*, \mu_t^*, \lambda_t^*)$. The SOSC is a regularity condition on the derivative of the KKT stationarity condition — the detailed definition of "SOSC" can be found in Appendix B. Then, one has the following proposition (Izmailov and Solodov, 2003; Izmailov et al., 2013).

**Proposition 2** (*upper Lipschitz stability of the solutions of KKT system under canonical perturbations*)

It is supposed that the SOSC is satisfied for the optimal solution $(x_t^*, \mu_t^*, \lambda_t^*)$ of the KKT system of *Problem* (1). Then, there exists a neighborhood $\mathcal{U}$ of $(x_t^*, \mu_t^*, \lambda_t^*)$ and $l > 0$, such that, for any $\sigma = (\sigma_1, \sigma_2, \sigma_3)$ close enough to $(0,0,0)$, any solution $(x_t(\sigma), \mu_t(\sigma), \lambda_t(\sigma)) \in \mathcal{U}$ of the perturbed KKT system

$$\nabla_{x_t} f_t(\hat{x}_t^*, \xi_t) + \nabla_{x_t} V_{t+1}(\hat{x}_t^*) + \sum_{i=1}^{p_t} \hat{\mu}_{t,i}^* \nabla_{x_t}\left(g_{t,i}(\hat{x}_t^*, \xi_t) + h_{t,i}(x_{t-1}, \xi_t)\right)$$

$$+ \sum_{j=1}^{q_t} \hat{\lambda}_{t,j}^* \nabla_{x_t}\left(l_{t,j}(\hat{x}_t^*, \xi_t) - b_{t,j}(x_{t-1}, \xi_t)\right) = \sigma_1,$$

$$g_t(x_t, \xi_t) + h_t(x_{t-1}, \xi_t) = \sigma_2$$

$$l_t(x_t, \xi_t) - b_t(x_{t-1}, \xi_t) = \sigma_3$$

satisfies the estimate

$$\|x_t(\sigma) - x_t^*\| + dist\left((\mu_t(\sigma), \lambda_t(\sigma)), (\lambda_t^*, \mu_t^*)\right) \leq l\|\sigma\|$$

**Proof.** See Property 1 of Izmailov et al. (2013).

Proposition 2 leads to the following theorem.

**Theorem 1** (*optimal convergence of solutions from the approximated value function*)
It is supposed that the SOSC is satisfied for the true optimal solution $(x_t^*, \mu_t^*, \lambda_t^*)$ of the KKT system. Let $\epsilon > 0$ be close enough to $0 \in R^{n_t}$ and with an associated $\mathcal{U}, l > 0$ from Proposition 2. Then, for sufficiently small $\epsilon_1 > 0$ and a solution $(\hat{x}_t^*, \hat{\mu}_t^*, \hat{\lambda}_t^*,) \in \mathcal{U}$ from the approximated value function, $\exists \epsilon_2 > 0$ such that

$$D_t(\theta_{t+1}; x_t) < \epsilon_2 \Rightarrow \|\hat{x}_t^* - x_t^*\| + dist\left((\hat{\mu}_t^*, \hat{\lambda}_t^*), (\mu_t^*, \lambda_t^*)\right) < \epsilon_1$$

**Proof.** See Appendix C.

Under the specified regularity condition in Theorem 1, the approximated solution $\hat{x}_t^*$ is guaranteed to converge to the true optimal solution $x_t^*$ as the stagewise KKT perturbation $D_t(\theta_{t+1}; x_t)$ approaches zero. Hence, VFGL should satisfy an additional condition (SOSC) for its theoretical convergence to the optimal solution, compared to those of well-known decomposition algorithms such as progressive hedging and SDDP. Because the MSLP has piecewise linear value functions that are not differentiable, they do not satisfy the SOSC. However, in Section 3, we show that the VFGL algorithm demonstrates a remarkable performance on MSLP problems even without the theoretical guarantee of convergence to the optimal solution.

Minimizing $D_t(\theta_{t+1}; x_t)$ is equivalent to minimizing its squared term $\{D_t(\theta_{t+1}; x_t)\}^2$. For computational reasons, the loss function $J_{t+1}(\theta_{t+1}; x_t)$ to minimize the estimation error $D_t(\theta_{t+1}; x_t)$ is defined as

$$
\begin{aligned}
J_{t+1}(\theta_{t+1}; x_t) &= \left\| \nabla V_{t+1}(x_t) - \nabla \hat{V}_{t+1}(x_t; \theta_{t+1}) \right\|^2 \\
&= \sum_{n=1}^{n_t} \left( \frac{\partial}{\partial x_{t,n}} V_{t+1}(x_t) - \frac{\partial}{\partial x_{t,n}} \hat{V}_{t+1}(x_t; \theta_{t+1}) \right)^2 \\
&= \sum_{n=1}^{n_t} \left( \frac{\partial}{\partial x_{t,n}} \mathbb{E}_{\xi_{t+1}}[\mathcal{V}_{t+1}(x_t, \xi_{t+1})] - \frac{\partial}{\partial x_{t,n}} \hat{V}_{t+1}(x_t; \theta_{t+1}) \right)^2
\end{aligned}
$$

The notation $\mathbb{E}_{\xi_{t+1}}[\cdot]$ is used to emphasize that the expectation operator is taken with respect to $\xi_{t+1}$.

The Leibnitz integral rule is applied to change the order of the partial derivative and expectation as follows.

$$J_{t+1}(\theta_{t+1}; x_t) = \sum_{n=1}^{n_t} \left( \mathbb{E}_{\xi_{t+1}} \left[ \frac{\partial}{\partial x_{t,n}} \mathcal{V}_{t+1}(x_t, \xi_{t+1}) \right] - \frac{\partial}{\partial x_{t,n}} \hat{V}_{t+1}(x_t; \theta_{t+1}) \right)^2$$

which is equivalent to

$$\sum_{n=1}^{n_t} \mathbb{E}_{\xi_{t+1}}\left[\left(r_{t+1,n}(\theta_{t+1}; x_t, \xi_{t+1})\right)^2\right] - Var_{\xi_{t+1}}(r_{t+1,n}(\theta_{t+1}; x_t, \xi_{t+1})) \tag{4}$$

where $r_{t+1,n}(\theta_{t+1}; x_t, \xi_{t+1}) = \frac{\partial}{\partial x_{t,n}} V_{t+1}(x_t, \xi_{t+1}) - \frac{\partial}{\partial x_{t,n}} \hat{V}_{t+1}(x_t; \theta_{t+1})$ is the residual between the true and estimated gradients, and $Var$ is the variance operator.

With elementary statistics, one can estimate Eq. (4) with $S$ samples $\xi_t^s, s = 1, \dots, S$.

$$\sum_{n=1}^{n_t} \left(\frac{1}{S}\sum_{s=1}^{S}\left(r_{t+1,n}(\theta_{t+1}; x_t, \xi_{t+1}^s)\right)^2 - \frac{1}{S-1}\sum_{s=1}^{S}\left(r_{t+1,n}(\theta_{t+1}; x_t, \xi_{t+1}^s) - \bar{r}_{t+1,n}(\theta_{t+1}; x_t)\right)^2\right)$$
$$\approx \frac{1}{S}\sum_{s=1}^{S}\left(\sum_{n=1}^{n_t} r_{t+1,n}(\theta_{t+1}; x_t, \xi_{t+1}^s)^2 - \left(r_{t+1,n}(\theta_{t+1}; x_t, \xi_{t+1}^s) - \bar{r}_{t+1,n}(\theta_{t+1}; x_t)\right)^2\right) \tag{5}$$

where $\frac{1}{S-1} \approx \frac{1}{S}$ is used for sufficiently large $S$, and

$\bar{r}_{t+1,n}(\theta_{t+1}; x_t) = \frac{1}{S}\sum_{s=1}^{S} r_{t+1,n}(\theta_{t+1}; x_t, \xi_{t+1}^s)$ is denoted as the sample mean of $r_{t+1,n}(\theta_{t+1}; x_t, \xi_{t+1})$.

The last term in Eq. (5) is defined as the approximated loss function $\mathcal{J}_{t+1}(\theta_{t+1}; x_t)$ at $x_t$. Then,

$$\mathcal{J}_{t+1}(\theta_{t+1}; x_t) = \frac{1}{S}\sum_{s=1}^{S}\left(\sum_{n=1}^{n_t} 2\bar{r}_{t+1,n}(\theta_{t+1}; x_t)r_{t+1,n}(\theta_{t+1}; x_t, \xi_{t+1}^s) - \bar{r}_{t+1,n}(\theta_{t+1}; x_t)^2\right)$$
$$= \frac{1}{S}\sum_{s=1}^{S} \jmath_{t+1}(\theta_{t+1}; x_t, \xi_t^s)$$

where $\jmath_{t+1}(\theta_{t+1}; x_t, \xi_t^s) = \sum_{n=1}^{n_t} 2\bar{r}_{t+1,n}(\theta_{t+1}; x_t)r_{t+1,n}(\theta_{t+1}; x_t, \xi_{t+1}^s) - \bar{r}_{t+1,n}(\theta_{t+1}; x_t)^2$

Note that

$$\nabla_{\theta_{t+1}} \bar{r}_{t+1,n}(\theta_{t+1}; x_t) = \frac{1}{S}\sum_{s=1}^{S} \nabla_{\theta_{t+1}} r_{t+1,n}(\theta_{t+1}; x_t, \xi_{t+1}^s)$$
$$= \frac{1}{S}\sum_{s=1}^{S} \nabla_{\theta_{t+1}}(\frac{\partial}{\partial x_{t,n}} V_{t+1}(x_t, \xi_{t+1}^s) - \frac{\partial}{\partial x_{t,n}} \hat{V}_{t+1}(x_t; \theta_{t+1}))$$
$$= -\frac{1}{S}\sum_{s=1}^{S} \nabla_{\theta_{t+1}} \frac{\partial}{\partial x_{t,n}} \hat{V}_{t+1}(x_t; \theta_{t+1})$$
$$= -\nabla_{\theta_{t+1}} \frac{\partial}{\partial x_{t,n}} \hat{V}_{t+1}(x_t; \theta_{t+1})$$

which is identical to the theta gradient on $r_{t+1,n}$ as follows.

$$\nabla_{\theta_{t+1}} r_{t+1,n}(\theta_{t+1}; x_t, \xi_{t+1}^s) = -\nabla_{\theta_{t+1}} \frac{\partial}{\partial x_{t,n}} \hat{V}_{t+1}(x_t; \theta_{t+1})$$

Therefore, the gradient of $\jmath_{t+1}(\theta_{t+1}; x_t, \xi_{t+1}^s)$ with respect to $\theta_{t+1}$ is given by

$$\nabla_{\theta_{t+1}} j_{t+1}(\theta_{t+1}; x_t, \xi_{t+1}^s) = -2 \sum_{n=1}^{n_t} r_{t+1,n}(\theta_{t+1}; x_t, \xi_{t+1}^s) \nabla_{\theta_{t+1}} \frac{\partial}{\partial x_{t,n}} \hat{V}_{t+1}(x_t; \theta_{t+1})$$

Because $x_t$ is adapted to history $\xi_{[t]}$, the stage $t+1$ loss function is defined by the expectation of $J_{t+1}(\theta_{t+1}; x_t)$ with respect to $x_t$ as follows.

$$\begin{aligned} J_{t+1}(\theta_{t+1}) &= \mathbb{E}_{x_t}[J_{t+1}(\theta_{t+1}; x_t)] \\ &\approx \frac{1}{N} \sum_{i=1}^{N} J_{t+1}(\theta_{t+1}; x_t^i) \end{aligned} \quad (6)$$

where $N$ is the number of samples of $x_t$, and $x_t^i$ is the $i$th sample of $x_t$.

The loss function in the form of Eq. (6) makes it possible to find a theta that minimizes the loss online based on the stochastic gradient descent (SGD). The application of the SGD is discussed in Section 3.3.

### 3.2.1. Sampling gradient of loss function

In this section, the process of sampling $J_{t+1}(\theta_{t+1}; x_t^i)$ is described. Because $\frac{\partial}{\partial x_{t,n}} \hat{V}_{t+1}(x_t; \theta_{t+1})$ and $\nabla_{\theta_{t+1}} \frac{\partial}{\partial x_{t,n}} \hat{V}_{t+1}(x_t; \theta_{t+1})$ can be easily computed, it remains to sample the target gradient $\nabla_{x_t} V_{t+1}(x_t, \xi_{t+1}^s)$. For a given pair of $(x_t, \xi_{t+1}^s)$, $\nabla_{x_t} V_{t+1}(x_t, \xi_{t+1}^s)$ can be found by differentiating the KKT optimality conditions of simulated subproblems (Kuhn & Tucker, 1950; Barratt, 2018). It is supposed that there is a primal dual optimal solution $(x_{t+1}^*, \lambda_{t+1}^*, \mu_{t+1}^*)$ for the stage $t+1$ subproblem defined in *Problem* (1) at the prior stage decision variable $x_t$ and realization $\xi_{t+1}^s$. Then, one can readily check from the Lagrangian

$$\begin{aligned} L(x_{t+1}^*, \lambda_{t+1}^*, \mu_{t+1}^* | x_t, \xi_{t+1}^s) &= f_t(x_{t+1}^*, \xi_{t+1}^s) + V_{t+2}(x_{t+1}^*) \\ &+ \sum_{i=1}^{p_{t+1}} \mu_{t+1,i}^* \left( g_{t+1,i}(x_{t+1}^*, \xi_{t+1}^s) + h_{t+1,i}(x_t, \xi_{t+1}^s) \right) \\ &+ \sum_{j=1}^{q_{t+1}} \lambda_{t+1,j}^* \left( l_{t+1,j}(x_{t+1}^*, \xi_{t+1}^s) + b_{t+1,j}(x_t, \xi_{t+1}^s) \right) \end{aligned} \quad (7)$$

that the local sensitivity of the subproblem objective with respect to $x_t$, at $\xi_{t+1}^s$ is the following:

$$\nabla_{x_t} V_{t+1}(x_t, \xi_{t+1}^s) = \sum_{i=1}^{p_{t+1}} \mu_{t+1,i}^* \nabla_{x_t} h_{t,i}(x_t, \xi_{t+1}^s) + \sum_{j=1}^{q_{t+1}} \lambda_{t+1,j}^* \nabla_{x_t} b_{t,j}(x_t, \xi_{t+1}^s) \quad (8)$$

Unfortunately, Eq. (8) is not available because finding the optimal solution of Eq. (7) requires the next stage value function $V_{t+2}$. Therefore, Eq. (8) is approximated with bootstrapping by replacing $V_{t+2}(x_{t+1})$ with the current best approximation of it, $\hat{V}_{t+2}(x_{t+1}; \theta_{t+1})$, which is analogous to the Benders cut approximation of SDDP.

### 3.2.2. Objective weighting

The loss function of VFGL is defined as a norm of multidimensional errors, which implies that the error from each dimension is treated equally in magnitude. However, each decision variable might have a different scale of sensitivity to the value function. For example, one can consider the asset liability management of public pension plans, minimizing the contribution rate while maximizing the wealth of the terminal stage. Typically, pension plan problems include decision variables of the contribution rate and allocation of wealth in various asset classes. A marginal increment in the contribution rate leads to an additional employee's total wage, whereas a similar nominal increment of allocation on Treasury securities affects the growth of wealth based on the risk-free rate. Apparently, the sensitivity of the value function with respect to Treasury securities allocation would be much lower than that of the contribution rate. A loss function without consideration of a decision variable sensitivity scale might lead to poor performance because decision variables with low sensitivity can be easily neglected in learning. In this regard, a weighted loss function is defined as

$$\mathcal{J}_{t+1}^w(\theta_{t+1}; x_t) = \sum_{n=1}^{n_t} \left( \frac{\frac{\partial}{\partial x_{t,n}} V_{t+1}(x_t) - \frac{\partial}{\partial x_{t,n}} \hat{V}_{t+1}(x_t; \theta_{t+1})}{w_{t+1,n}} \right)^2,$$

where $w_{t+1,n} = \mathbb{E}_{x_t}\left[\frac{\partial}{\partial x_{t,n}} V_{t+1}(x_t)\right], n = 1, \ldots, n_t$.

Subsequently, one obtains a scaled loss gradient as follows.

$$\nabla_{\theta_{t+1}} \hat{\jmath}_{t+1}^w(\theta_{t+1}; x_t, \xi_{t+1}^s) = -2 \sum_{n=1}^{n_t} \frac{1}{w_{t+1,n}^2} \left( \frac{\partial}{\partial x_{t,n}} V_{t+1}(x_t, \xi_{t+1}^s) - \frac{\partial}{\partial x_{t,n}} \hat{V}_{t+1}(x_t; \theta_{t+1}) \right) \nabla_\theta \frac{\partial}{\partial x_{t,n}} \hat{V}_{t+1}(x_t; \theta_{t+1})$$

However, the value of $w_{t+1,n}$ cannot be computed directly. Therefore, the value $w_{t+1,n}$ is replaced with its sample mean $\bar{w}_{t+1,n}$, and this approximation is improved as the iteration continues — that is, whenever $\frac{\partial}{\partial x_{t,n}} \hat{V}_{t+1}(x_t; \theta_{t+1})$ is sampled at iteration $i$, the following is updated:

$$\bar{w}_{t+1,n} \leftarrow \frac{i}{i+1} \bar{w}_{t+1,n} + \frac{1}{i+1} \left( \frac{\partial}{\partial x_{t,n}} \hat{V}_{t+1}(x_t; \theta_{t+1}) \right)$$

Because the major purpose of $\bar{w}_{t+1,n}$ is scaling, an inaccurate approximation of $w_{t+1,n}$ is not a serious concern as long as it reflects the average scale.

### 3.3. Parameter optimization

Parameters $\theta_{t+1}$ that minimize the weighted loss functions $J_{t+1}^w(\theta_{t+1})$ for $t = 1, \ldots T - 1$ are optimized with the SGD algorithm. The SGD algorithm for VFGL is described in Algorithm 1.

---

**Algorithm 1:** Stochastic Gradient Descent for VFGL

**Require:** $N$ (maximum iteration), $\alpha_i$ for $i = 1, \ldots, N$ (step size), $S$ (sample size)
**Initialize:** $i \leftarrow 1$ (iteration count), $\theta_t\ for\ t = 1, \ldots, T$ (parameter values)

1   while $i \leq N$ do
2      for $t = 1, \ldots, T$ do
3         if $t = 1$ do
4            Solve the first stage subproblem of *Problem* (2). Save optimal solution $x_1^i$.
5         else if $t > 1$ do

```
6        for s = 1,…S do
7            Sample ξₜˢ
8            Solve the stage t subproblem of Problem (2). Save optimal solution (xₜˢ, λₜˢ, μₜˢ).
9            Compute ∇_{θₜ} 𝒥ₜʷ(θₜ; xₜ₋₁ⁱ, ξₜˢ)
10       end for
11       Compute ∇_{θₜ}𝒥ₜʷ(θₜ; xₜ₋₁ⁱ) = (1/S)∑_{s=1}^{S} ∇_{θₜ}𝒥ₜʷ(θₜ; xₜ₋₁ⁱ, ξₜˢ)
12       Update θₜ ← θₜ − αᵢ∇_{θₜ}𝒥ₜʷ(θₜ; xₜ₋₁ⁱ)
13       xₜⁱ ← randomly chosen from {xₜ¹, …, xₜᵏ}
14    end if
15  end for
16  i ← i + 1
17 end while
```

Here, the step size (or learning rate) $\alpha_i$ satisfies

$$\sum_{i=1}^{\infty} \alpha_i = \infty, \sum_{i=1}^{\infty} \alpha_i^2 < \infty$$

for the convergence (Robbins & Monro, 1951).

SGD is widely used in the parameter optimization of artificial intelligence algorithms, and there are a number of variants of SGD to address problem-specific issues. One variant of SGD that can be particularly useful for this problem is implicit stochastic gradient descent (ISGD). ISGD involves the extra step of calculating the updated parameter implicitly defined by the following equation.

$$\theta_t^{new} = \theta_t^{old} - \alpha_i \nabla_{\theta_t} \mathcal{J}_t^w(\theta_t^{new}; x_{t-1})$$

Despite the additional computational burden, the ISGD is known to have a stable convergence regardless of the choice of learning rate $\alpha_i$. The ISGD is suitable for VFGL because the unstable scale of gradients of VFGL tends to require a careful selection of the learning rate for convergence, which is cumbersome.

In Algorithm 1, $\nabla_{\theta_t}\mathcal{J}_t^w(\theta_t; x_{t-1}^i)$ is approximated with $S$ samples. One can employ some variance reduction techniques (Hammersley & Morton, 1956; Siegmund, 1976; Rubinstein & Marcus, 1985; Haghighat & Wagner, 2003) to enhance the sample approximation. In particular, moment matching is used. The moment matching algorithms vary the values of finite samples and/or probability assigned to them so that the moments of the samples are matched to that of their original distribution.

The following procedure can be considered one of the most naïve forms of the moment matching algorithm. For $S$ samples $y_1, …, y_S$ drawn from their true distribution $Y$ with mean $\mu$ and covariance $\Sigma$, each $y_s$ is transformed as follows.

$$z_s = \mu + \Sigma^{\frac{1}{2}} C^{-\frac{1}{2}}(y_s - \bar{y}),$$

where $\bar{y}$ is the sample mean, and $C$ is the sample covariance matrix. One can easily check that the mean and variance of $z_s$ exactly match those of $Y$. However, matching up to the second moment might not be sufficient, depending on the underlying distribution and structure of the problems. For more details of moment matching algorithms, see Høyland et al. (2003), Ji et al. (2005), Ponomareva et al. (2015), Staino & Russo (2015), and Radhakrishnan et al. (2018).

### 3.3.1. Stopping criterion

The convergence of parameter values can be used as the stopping criterion of VFGL. One cannot employ the expected optimality gap for the stopping criterion as in SDDP because a first stage problem in VFGL no longer provides a lower bound of the problem. However, one can observe the change in parameters of the value function. In particular, a total parameter change for the $i$-th iteration is defined as follows.

$$\Delta \theta^i = \sum_{t=2}^{T} \Delta \theta_t^i,$$

where $\Delta \theta_t^i = \left\| \theta_t^i - \theta_t^{i-1} \right\|_2$ is the parameter change for the stage $t$ value function $V_t$ in the $i$-th iteration. When $\Delta \theta^i$ is sufficiently small, one can conclude that the value function converges. The entire procedure of VFGL is presented in Algorithm 2.

---

**Algorithm 2**: VFGL

**Require:** $N$ (maximum iteration), $\alpha_i$ for $i = 1, ..., N$ (step size),
$S$ (sample size), $\epsilon$ (minimum theta update threshold),
$\hat{V}_{t+1}$ for $t = 1, ..., T - 1$ (parametric value function),

**Initialize:** $i \leftarrow 1$ (iteration counter)
$\Delta \theta \leftarrow \infty$ (total theta update level)
$p_t^s \leftarrow \frac{1}{S}$
$\theta_t^1$ for $t = 2, ..., T - 1$ (parameter values)

1  while $\Delta \theta^i > \epsilon$ and $i < N$ do
2     $\Delta \theta^i \leftarrow 0$
3     for $t = 1, ..., T$ do
4       if $t = 1$ do
5         Solve the first stage subproblem of *Problem* (2).
6         Save optimal solution $x_1^i$.
7       else if $t > 1$ do
8         Sample random variables $\xi_t^s$ for $s = 1, ..., S$
9         *(optional)* Match moments of sampled $\xi_t^s$ with the corresponding probability $p_t^s$
10        for $s = 1, ..., S$ do
11          Solve the corresponding subproblem of *Problem* (2)
12          Compute the loss gradient $\nabla_{\theta_t} \mathcal{J}_t^w(\theta_t^i; x_{t-1}^i, \xi_t^s)$
13        Compute the minibatch loss gradient $\nabla_{\theta_t} \mathcal{J}_t^w(\theta_t^i; x_{t-1}^i) = \sum_{s=1}^{S} p_t^s \nabla_{\theta_t} \mathcal{J}_t^w(\theta_t^i; x_{t-1}^i, \xi_t^s)$
14        Update $\theta_{t+1}$ by ISGD on $\theta_{t+1}^{i+1} = \theta_{t+1}^i - \alpha_i \nabla_{\theta_t} \mathcal{J}_t^w(\theta_t^i; x_{t-1}^i)$
         Sample $\xi_t^i$ solve the stage $t$ subproblem of *Problem* (2) and save $x_t^i$
15        $\Delta \theta^i \leftarrow \Delta \theta^i + \left\| \theta_t^{i+1} - \theta_t^i \right\|_2$
16      end if
17    end for
18    $i \leftarrow i + 1$

### 3.4. Quality measure of parametric form: KKT deviation

The biggest hurdle of the proposed algorithm is finding an appropriate parametric form of the value function. As pointed out by Agrawal et al. (2020), the tractable form of the true value function can be found only in a few cases. Fortunately, the choice of the parametric approximation of the value function $\hat{V}_t$ need not be exact. If the modeling capacity of $\hat{V}_t$ is strong enough that its gradients can be sufficiently close to the gradients of the true value, VFGL can find a solution that is sufficiently close to the true optimal solution. There may be multiple parametric forms of the value function that can well approximate the true value function, and they must be compared.

In this regard, a measure is proposed to evaluate the performance of various parametric approximations to choose the best one among them. The solution of VFGL converges to the true solution as $D_t$ goes to 0 by Theorem 1. Let $D$ be the expected value of the sum of $D_t$ over stages as follows.

$$D(\theta_2, \ldots, \theta_T) = \mathbb{E}\left[\sum_{t=1}^{T-1} D_t(\theta_{t+1}; x_t)\right]$$

One forward simulates solution $x_t^i, t = 1, \ldots, T \ \ i = 1, \ldots, N$ based on the current approximation $\hat{V}_{t+1}, t = 1, \ldots T-1$ to obtain the estimator of $D$ as follows.

$$\hat{D}(\theta_2, \ldots, \theta_T) = \frac{1}{N}\sum_{i=1}^{N}\sum_{t=1}^{T-1}\left\|\frac{1}{S}\sum_{s=1}^{S}\nabla_{x_t}\mathcal{V}_{t+1}(x_t^i, \xi_{t+1}^s) - \nabla_{x_t}\hat{V}_{t+1}(x_t^i; \theta_{t+1})\right\|_2$$

where each $\nabla_{x_t}\mathcal{V}_{t+1}(x_t^i, \xi_{t+1}^s)$ is sampled from the bootstrapping described in Section 2.2.1.

Here, $\hat{D}$ is referred to as KKT deviation, and it indicates the performance of the value function gradient approximation based on the current parametric form. A better quality of the solution would be expected for a low KKT deviation, but a high KKT deviation does not necessarily lead to an erroneous solution. Instead, the KKT deviation can be regarded as a confidence level of the solution. Illustrative examples of the relationship between the KKT deviations and solution qualities are presented in Section 4.

**Remark 1.** A good candidate for the initial choice of the parametric form of the value function is the functional form of the stagewise objective function or the indefinite integral form of the stagewise objective function. Such a claim is based on the observation that the past stage decision variables of the decomposed multistage stochastic programming problems are frequently used as a resource constraint. Furthermore, the stagewise objective function is frequently a function of such available resources. However, sometimes the analytical form of the indefinite integral of the objective function might not be available. In such a case, one can approximate the indefinite integral form with sampling. More detailed illustrations and discussions can be found in Section 4.

## 4. Illustrative examples

In this section, three illustrative examples with different application domains are presented to demonstrate the empirical performance of VFGL. The three examples are financial planning, production planning, and hydrothermal energy planning, where they have different types of decision variables, objectives, and constraints. These problems have been extensively studied within the optimization community.

In the numerical experiments, seven to eleven stage multistage stochastic convex problems are considered to compare the results from MSP, SDDP, and VFGL. Here, MSP is a deterministic equivalent of the multistage stochastic convex problems. MSP was used as a benchmark to assess the solution quality of SDDP and VFGL because MSP finds a global optimum given a finite scenario tree. The stopping criteria for SDDP and VFGL are different. Therefore, the two algorithms were run until solution convergence, and their first stage solution, elapsed time, and degree of convergence at each iteration were compared. For SDDP, 20 scenarios for upper bound calculation in the forward pass were sampled, but the backward pass was performed on only one sampled path. All the experiments were performed using an Intel i5-7500 processor with 32 GB of RAM. For each stagewise subproblem, MOSEK[3] solver 12.9.0 was used to solve the linear and convex problems by means of CVXPY 1.0.28, a Python-based modeling framework (Diamond & Boyd, 2016).

### 4.1 Lifetime financial planning

Merton (1969) solved a continuous-time portfolio optimization problem, where the objective is to maximize the total utility of consumption and terminal wealth by appropriately deciding the level of consumption and the allocation of wealth in stocks and bonds. See Appendix E for the full problem definition and its analytical solution. For the numerical experiment, the time period is discretized into $T$ stages, and the problem is approximated using an MSP formulation. The approximated problem is solved with MSP, SDDP, and VFGL. The variables, parameters, and parameter value for the discretized lifetime financial planning problem case study are described in Table 1.

*Table 1. Decision variables and parameters for discretized lifetime financial planning*

| Decision variables | Description | |
|---|---|---|
| $C_t$ | Consumption at stage $t$ | |
| $W_t$ | Wealth in the beginning of stage $t$ | |
| $S_t$ | Amount invested into stock at stage $t$ | |
| $B_t$ | Amount invested into bond at stage $t$ | |

| Parameters | Description | Value |
|---|---|---|
| $\rho$ | Discount rate. | 0 |

---

[3] See http://www.mosek.com

| $\gamma$ | Utility risk aversion coefficient | 1 |
| $\mu, \sigma$ | Mean return, volatility respectively, of stock | 1.06, 0.20 |
| $r$ | Risk-free rate of bond | 1.03 |
| $\epsilon$ | Scaling coefficient of bequest utility | 1 |
| $\xi_t$ | Random variable following a standard normal distribution | - |

The total cost is as follows.

$$-\mathbb{E}\left[\sum_{t=1}^{T} e^{-\rho t} U(C_t) + \epsilon^\gamma e^{-\rho T} U(W_T)\right]$$

Stagewise subproblems are defined as follows.

***Stage 1 subproblem***

$$\begin{aligned}
\text{minimize} \quad & -U(C_1) + V_2(S_1, B_1) \\
\text{subject to} \quad & S_1 + B_1 + C_1 = 1 \quad &&\text{Initial wealth} \\
& S_1, B_1, C_1 \geq 0 \quad &&\text{Non-negativity}
\end{aligned}$$

***Stage t subproblem*** $(t = 2, \ldots, T-1)$

$$\begin{aligned}
\text{minimize} \quad & -U(C_t) + V_{t+1}(S_t, B_t) \\
\text{subject to} \quad & S_t + B_t + C_t = r\Delta t B_{t-1} + e^{\left(\mu - \frac{\sigma^2}{2}\right)\Delta t + \sigma\sqrt{\Delta t}*\xi_t} S_{t-1} \quad &&\text{Inventory balance} \\
& S_t, B_t, C_t \geq 0 \quad &&\text{Non-negativity}
\end{aligned}$$

***Stage T subproblem***

$$\begin{aligned}
\text{minimize} \quad & -U(W_T) \\
\text{subject to} \quad & W_T = r\Delta t B_{T-1} + e^{\left(\mu - \frac{\sigma^2}{2}\right)\Delta t + \sigma\sqrt{\Delta t}*\xi_T} S_{T-1} \quad &&\text{Inventory balance}
\end{aligned}$$

From the special structure of the terminal stage problem, $V_T(S_{T-1}, B_{T-1})$ can be calculated as follows.

$$V_T(S_{T-1}, B_{T-1}) = \int_{\xi_T} -\ln\left(r\Delta t B_{T-1} + e^{\left(\mu - \frac{\sigma^2}{2}\right)dt + \sigma\sqrt{\Delta t}*\xi_T} S_{T-1}\right) \varphi(\xi_T) d\xi_T$$

where $\varphi(\xi_T)$ is a probability density function of $\xi_T$. Finding an exact parametric form of $V_T(S_{T-1}, B_{T-1})$ above is extremely difficult. Instead, the parametric form of $V_T(S_{T-1}, B_{T-1})$ is approximated with finite sampling as follows.

$$\int_{\xi_T} -\ln\left(rdt B_{T-1} + e^{\left(\mu - \frac{\sigma^2}{2}\right)dt + \sigma\sqrt{dt}*\xi_T} S_{T-1}\right) \varphi(\xi_T) d\xi_T \approx \frac{1}{b_T} \sum_{i=1}^{b_T} -ln\left(rB_{T-1} + \beta_{i,T} S_{T-1}\right)$$

where $\beta_{i,T}$ is a sampled value of $e^{\left(\mu - \frac{\sigma^2}{2}\right)\Delta t + \sigma\sqrt{\Delta t}*\xi_T}$, and $b_T$ is the sampling number. Therefore, the parametric value function is set to include the parametric term as follows.

$$\theta_t^1 \sum_{i=1}^{b_T} -\ln(rB_{T-1} + \beta_{i,T} S_{T-1})$$

A negative utility term is added for flexibility by $-\theta_t^2 \ln(\theta_t^3 B_t + \theta_t^4 S_t)$ to increase the model capacity. Therefore, the parametric value function for stage $t = 2, \dots, T$ is defined as follows.

$$\hat{V}_t(S_{t-1}, B_{t-1}; \theta_t = (\theta_t^1, \theta_t^2, \theta_t^3, \theta_t^4)) = -\theta_t^1 \sum_{i=1}^{b_t} \ln(rB_{t-1} + \beta_{i,t} S_{t-1}) - \theta_t^2 \ln(\theta_t^3 B_{t-1} + \theta_t^4 S_{t-1})$$

where $b_t$ is the sampling number. Moreover, $b_t = 30$ is set for all t, and the parameters are initialized as follows.

$$\theta_t^1, \theta_t^2, \theta_t^3, \theta_t^4 = \left(\frac{1}{b_t}, \frac{1}{b_t}, 1, 1\right)$$

An 11-stage scenario tree is constructed using three samples per stage, where the samples are normalized to match the first and second moments of the true distribution perfectly. The constructed scenario tree is shared between MSP and SDDP. The experimental results for MSP, SDDP, and VFGL are summarized in Table 2. It presents the averaged first stage solutions over 20 independently generated instances of MSP, SDDP, and VFGL. The results from Table 2 verifies that the MSP solution matches well with the analytical solution, and VFGL provides more accurate solutions than SDDP.

*Table 2. Comparison of objective value, first stage solution, and computation time between algorithms for lifetime financial planning problem*

| Algorithm | Objective value | Risky asset portion | Risk-free asset | Risky asset | Consumption | Average computation time (s) |
|---|---|---|---|---|---|---|
| **Analytical** | - | 0.7500 | 0.2292 | 0.6875 | 0.0833 | - |
| **MSP** | 29.5508 (0.0000) | 0.7510 (0.000) | 0.2282 (0.000) | 0.6885 (0.000) | 0.0833 (0.000) | - |
| **SDDP** | 29.4377 (0.0577) | 0.7694 (0.012) | 0.2114 (0.012) | 0.7053 (0.012) | 0.0833 (0.000) | 4552 (15) |
| **VFGL** | 29.5112 (0.0305) | 0.7508 (0.000) | 0.2284 (0.000) | 0.6883 (0.000) | 0.0833 (0.000) | 354 (1) |

The evolution of the first stage solutions of (a) SDDP and (b) VFGL for a randomly generated instance is depicted in Figure 1. While the first stage solution of VFGL quickly converges, almost within 20 iterations, the solution of SDDP is not stabilized within 200 iterations. SDDP was further run for 1000 iterations, which took approximately 28,000 s, but the first stage solution was still not fully stabilized. The SDDP algorithm shows a particularly slow convergence of this example.

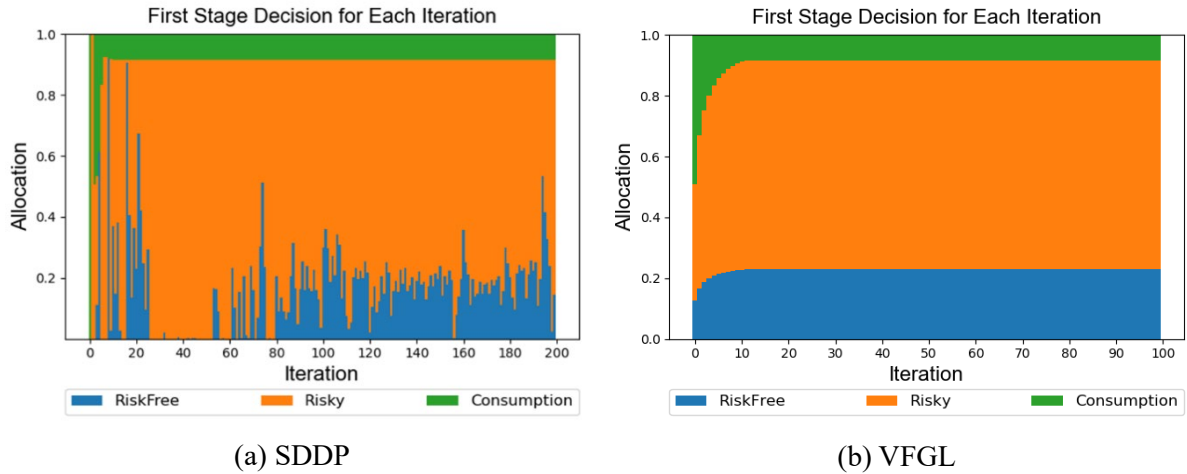

(a) SDDP  (b) VFGL

**Figure 1**. Comparison of evolution of the first stage solution between SDDP and VFGL for discretized lifetime financial planning problem

The computation time per iteration is illustrated in Figure 2. Figures 1 and 2 clearly show that VFGL is not only faster in convergence but also faster in each iteration than SDDP. Furthermore, SDDP shows increasing computation time with iteration. Because SDDP requires a high number of iterations for convergence, the computational advantage of VFGL over SDDP is significant in this sample.

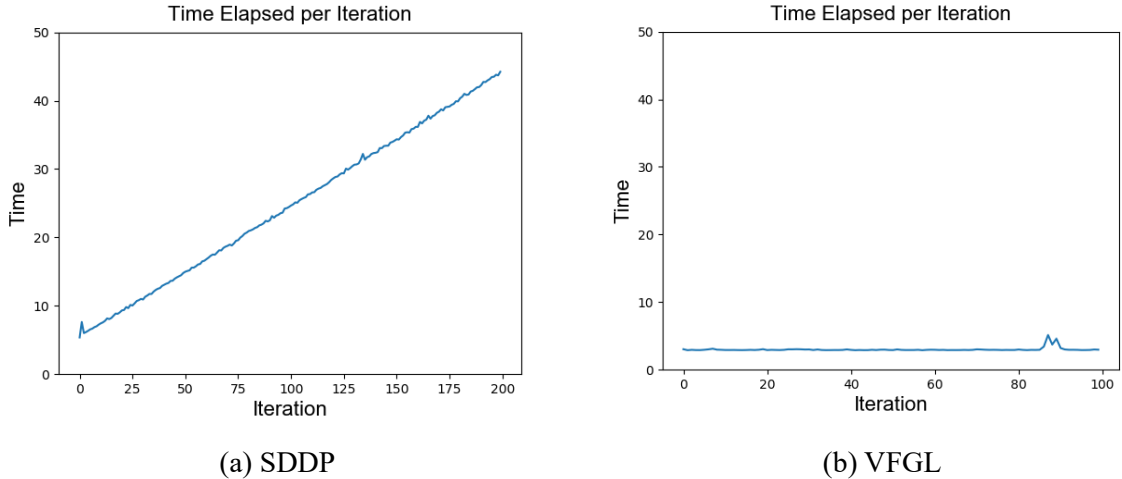

(a) SDDP  (b) VFGL

**Figure 2.** Comparison of computation time per iteration between SDDP and VFGL for discretized financial planning problem

The stopping criteria of SDDP and VFGL are displayed in Figure 3. The dashed upper bound of Figure 3(a) refers to the conservative upper bound estimation with 95% confidence (Shapiro, 2011). Figure 3(b) describes the total $\theta$ change in VFGL at each iteration. The convergence of the parameter for VFGL is apparent in Figure 3(b) after approximately 10 iterations, which coincides with the solution convergence shown in Figure 1. However, the approximated SDDP bound of Figure 3(a) provides little information about the solution convergence.

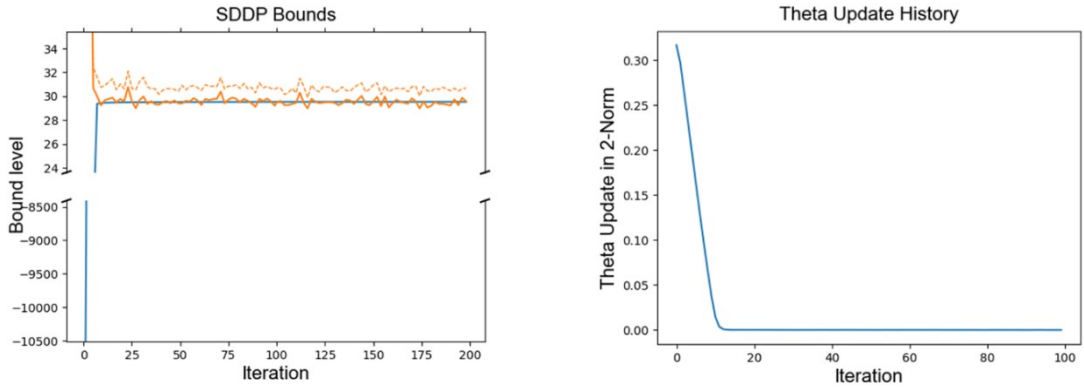

(a) SDDP bounds            (b) VFGL parameter update history

**Figure 3**. Comparison of stopping criteria between SDDP and VFGL for discretized lifetime financial planning problem

Different parametric forms of the value function are applied to this example, and the result changes are analyzed. By constructing the parametric form of the value function for this example, one can easily change the parametric family by changing the sampling number $b_t$. The KKT deviation, solution error, and objective value for different $b_t$ values are described in Table 3. The KKT deviation tends to decrease as $b_t$ increases. While a negative correlation between the KKT deviation and the first stage solution error is evident, a negative correlation between the KKT deviation and the objective value is not very clear in Table 3. An extremely flat objective function and objective approximation error may be the causes. Therefore, one can argue that the KKT deviation is a valid metric that indicates the possible appropriateness of the parametric form of the value function.

*Table 3. Comparison of KKT deviation, first stage solution error, and objective value between various sampling numbers for discretized financial planning problem (\* indicates the baseline case for the comparison)*

| Sampling number $b_t$ | KKT deviation | First stage solution error | Objective value |
| --- | --- | --- | --- |
| 2 | 6.7044 | 0.3258 | 29.6344 |
| 3 | 3.0451 | 0.3217 | 29.5369 |
| 4 | 1.0743 | 0.3155 | 29.5102 |
| 5 | 0.2846 | 0.2990 | 29.4814 |
| 10 | 0.0128 | 0.0399 | 29.5106 |
| 15 | 0.0103 | 0.0062 | 29.5396 |
| 20 | 0.0157 | 0.0013 | 29.4975 |
| 25 | 0.0159 | 0.0008 | 29.4246 |
| *30 | 0.0075 | 0.0001 | 29.5112 |

Finally, the risk-free return and standard deviation of the risky asset of the original problem are perturbed, and whether the given parametric form of the value function can be robustly recycled with respect to the perturbation is determined. The results of the perturbed problem are reported in Table 4. The table shows that both the solution and the objective from VFGL and MSP for any perturbation coincide within an approximation error. It can be concluded that the given parametric family of the value function can be recycled for a small perturbation of the asset distribution for this example.

*Table 4. Comparison of perturbed problems with a fixed parametric value function form for discretized lifetime financial planning problem (\* indicates the baseline case for the comparison)*

| Parameter | | | VFGL | | | | | MSP | | | |
|---|---|---|---|---|---|---|---|---|---|---|---|
| Risk-free return | Risky return | Risky volatility | KKT deviation | Objective | Risk free | Risky | Consumption | Objective | Risk free | Risky | Consumption |
| 0.04 | 0.060 | 0.20 | 0.0114 | 29.5281 | 0.4567 | 0.4600 | 0.0833 | 29.5266 | 0.4564 | 0.4603 | 0.0833 |
| 0.034 | 0.060 | 0.26 | 0.0105 | 29.5245 | 0.5122 | 0.4045 | 0.0833 | 29.5458 | 0.5120 | 0.4046 | 0.0833 |
| 0.031 | 0.060 | 0.29 | 0.0107 | 29.5399 | 0.5367 | 0.3799 | 0.0833 | 29.5565 | 0.5372 | 0.3795 | 0.0833 |
| *0.030 | *0.060 | *0.20 | 0.0075 | 29.5112 | 0.2284 | 0.6883 | 0.0833 | 29.5508 | 0.6885 | 0.2282 | 0.0833 |
| 0.028 | 0.060 | 0.32 | 0.0106 | 29.5429 | 0.5571 | 0.3596 | 0.0833 | 29.5667 | 0.5580 | 0.3587 | 0.0833 |
| 0.025 | 0.060 | 0.35 | 0.0108 | 29.5567 | 0.5746 | 0.3421 | 0.0833 | 29.5768 | 0.5755 | 0.3412 | 0.0833 |
| 0.022 | 0.060 | 0.38 | 0.0109 | 29.5703 | 0.5900 | 0.3267 | 0.0833 | 29.5867 | 0.5910 | 0.3257 | 0.0833 |
| 0.019 | 0.060 | 0.41 | 0.0110 | 29.5791 | 0.6036 | 0.3130 | 0.0833 | 29.5966 | 0.6047 | 0.3119 | 0.0833 |
| 0.016 | 0.060 | 0.44 | 0.0110 | 29.5905 | 0.6159 | 0.3008 | 0.0833 | 29.6065 | 0.6171 | 0.2996 | 0.0833 |
| 0.013 | 0.060 | 0.47 | 0.0111 | 29.5993 | 0.6270 | 0.2897 | 0.0833 | 29.6163 | 0.6282 | 0.2884 | 0.0833 |
| 0.01 | 0.060 | 0.5 | 0.0112 | 29.6081 | 0.6371 | 0.2795 | 0.0833 | 29.6261 | 0.6384 | 0.2783 | 0.0833 |

## 4.2 Production planning

A simple $T$-stage factory production/storage planning problem was considered. Specifically, a dynamic lot size problem was solved to find an optimal way of satisfying stochastic demand while minimizing the cost incurred by manufacturing, ordering, or carrying inventory. The problem is regarded as a standard production planning problem and has been widely studied (e.g., Wagner & Whitin, 1958; Shapiro, 1993; Karimi et al., 2003).

In the example, the factory can produce and store each product $i \in I$ while facing uncertain demand. In the first stage, there is no demand, and the factory manager must only decide how much to produce and store. It is assumed that the available resource for the production is fixed and cannot be transferred to the next stage. In addition, there is no ordering cost for production. In the second stage and forward, the manager observes realized demand and adjusts the production level to meet the demand. Any products left after demand are stored in the next stage with storage cost. When the demand exceeds the sum of supply and initial storage, the shortage is outsourced at a relatively high cost. An optimal production policy is sought that minimizes the sum of the expected storage and outsourcing costs. The example involves three products for seven stages. The demand is independent and identically distributed over stages with three possible values with equal probability. The decision variables and parameters, and values for parameters are listed in Table 5.

*Table 5. Decision variables and parameters for production planning problem*

| Decision variables | Description | |
|---|---|---|
| $x_{t,i}$ | Quantity of product $i$ produced at stage $t$ | |
| $y_{t,i}$ | Quantity of product $i$ outsourced at stage $t$ | |
| $s_{t,i}$ | Quantity of product $i$ stored at the end of stage $t$ | |
| **Parameters** | **Description** | **Value** |
| $a_{t,i}$ | Production cost of product $i$ at stage $t$ | (1, 2, 5) |
| $b_{t,i}$ | Outsourcing cost of product $i$ at stage $t$ | (6, 12, 20) |
| $c_{t,i}$ | Storage cost of product $i$ from the end of stage $t$ to beginning of stage $t+1$ | (3, 7, 10) |

| $r_t$ | Maximum production resource available at stage $t$ | 10 |
| $d_{t,i}$ | Random demand of product $i$ at stage $t$ | $\{(5,3,1),(6,2,1),(1,2,2)\}$ |

The total cost is as follows.

$$\mathbb{E}\left[\sum_{t=1}^{T-1}\sum_{i\in I}(y_{t,i}b_{t,i}+s_{t,i}c_{t,i})+y_{T,i}b_{T,i}\right]$$

The stagewise subproblems at the first stage, middle stages, and final stage are defined as follows.

***Stage 1 subproblem***

$\text{minimize}\quad \sum_{i\in I}y_{1,i}b_{1,i}+\sum_{i\in I}s_{1,i}c_{1,i}+v_2(s_1)$

$\text{subject to}\quad \sum_{i\in I}x_{1,i}a_{1,i}\leq r_1$      For $i\in I$      Resource limit

$\qquad\qquad\quad s_{1,i}=x_{1,i}+y_{1,i}$      For $i\in I$      Storage balance

$\qquad\qquad\quad x_{1,i},y_{1,i},s_{1,i}\geq 0$      For $i\in I$      Non-negativity

***Stage $t$ subproblem $(t=2,\ldots,T-1)$***

$\text{minimize}\quad \sum_{i\in I}y_{t,i}b_{t,i}+\sum_{i\in I}s_{t,i}c_{t,i}+v_{t+1}(s_t)$

$\text{subject to}\quad \sum_{i\in I}x_{t,i}a_{t,i}\leq r_t$      For $i\in I$      Resource limit

$\qquad\qquad\quad s_{t,i}=s_{t-1,i}+x_{t,i}+y_{t,i}-d_{t,i}$      For $i\in I$      Storage balance

$\qquad\qquad\quad x_{t,i},y_{t,i},s_{t,i}\geq 0$      For $i\in I$      Non-negativity

***Stage $T$ subproblem***

$\text{minimize}\quad \sum_{i\in I}y_{T,i}b_{T,i}$

$\text{subject to}\quad \sum_{i\in I}x_{T,i}a_{T,i}\leq r_t$      For $i\in I$      Resource limit

$\qquad\qquad\quad s_{T,i}=s_{T-1,i}+x_{T,i}+y_{T,i}-d_{T,i}$      For $i\in I$      Storage balance

$\qquad\qquad\quad x_{T,i},y_{T,i},s_{T,i}\geq 0$      For $i\in I$      Non-negativity

Note that the stagewise value functions of this problem are convex piecewise linear because each stagewise subproblem is linear programming. Since a piecewise linear function is not differentiable, the convergence of VFGL to the optimal is not guaranteed by the Theorem 1. Through this example, we limitedly check the viability of VFGL to multistage stochastic convex programming without the SOSC condition.

For all products, the outsourcing cost is significantly higher than the production resource cost and storage cost. Therefore, one obtains a crude insight that the optimal solution would avoid outsourcing products, which is likely to be achieved by creating nonempty storage. Furthermore, the storage cost is paid in advance. Therefore, one can expect that the value function would decrease with storage to some extent. However, the disutility of storage increases after a certain point because there exists a certain level of storage sufficient for providing a buffer for uncertain demands. From these points of view, one can use the following form of the parametric value function:

$$\hat{V}_t(s_{t-1}; \theta_t = (\theta_t^1, \theta_t^2)) = s_{t-1} \cdot \theta_t^1 + \sum_{i=1}^{3} e^{-\theta_{t,i}^2 s_{t-1,i}} \quad \text{for } t = 2, \ldots, T$$

Here, $\theta_t^1, \theta_t^2 \in \mathbf{R}^3$. The first linear term represents the overall decreasing disutility of storage, while the exponential term is added to reflect the increasing disutility of storage after a certain point. The parameters are initialized first using the following values.

$$\theta_t^1 = (-1, -1, -1), \theta_t^2 = (0, 0, 0) \text{ for } t = 2, \ldots, T$$

*Table 6. Comparison of objective value, first stage solution, and computation time between algorithms for production planning problem*

| Algorithm | Objective value | Production | | | Outsource | | | Computation time (s) |
|---|---|---|---|---|---|---|---|---|
| | | Product 1 | Product 2 | Product 3 | Product 1 | Product 2 | Product 3 | |
| Optimal (MSP) | 178 (0.0) | 1.00 (0.00) | 0.00 (0.00) | 1.00 (0.00) | 0.00 (0.00) | 0.00 (0.00) | 0.00 (0.00) | - |
| SDDP | 179 (0.4) | 1.00 (0.00) | 0.00 (0.00) | 1.02 (0.06) | 0.00 (0.00) | 0.00 (0.00) | 0.00 (0.00) | 1016 (5) |
| VFGL | 178 (0.2) | 0.99 (0.00) | 0.00 (0.00) | 1.01 (0.00) | 0.00 (0.00) | 0.00 (0.00) | 0.00 (0.00) | 157 (1) |

Table 6 presents the first stage solution of the MSP model and the averaged first stage solutions over 20 generated instances of SDDP and VFGL. SDDP took approximately 100 iterations for the first stage solution convergence, while VFGL took around 200 iterations. Parentheses in Table 6 indicate the standard deviation of each solution. Both SDDP and VFGL show accurate solutions and objective values within an acceptable standard deviation. However, there is a striking difference in the computation time between the two. The computation time for SDDP is seven times longer than that of VFGL, even with 100 fewer iterations. Figure 4 shows the first stage solution of (a) SDDP and (b) VFGL. Only production decision variables are shown because the first stage outsourcing decision variable always turned out to be 0, making the storage decision identical to the production decision. One can see that the first stage solution has been stabilized for both algorithms. Regarding the computation time and stopping criterion, similar findings to Section 3.1 can be seen in Appendix F.

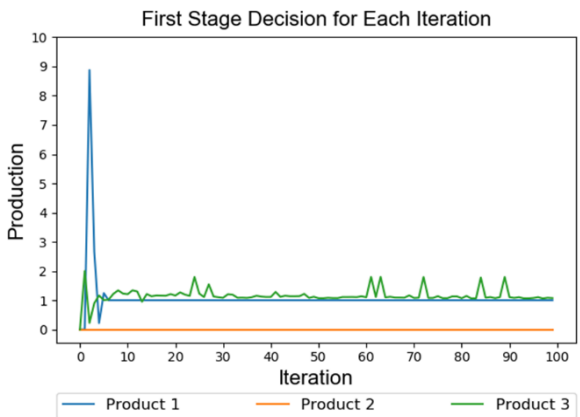 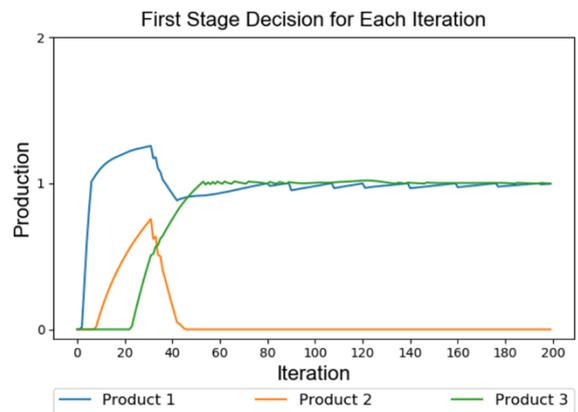

a) SDDP         (b) VFGL

**Figure 4**. Comparison of evolution of the first stage solution between SDDP and VFGL for production planning problem

To see the effect of the choice of parametric form in VFGL, a few more parametric forms were tested. The KKT deviation and first stage solution error for each parametric form of the value function are described in Table 7. The KKT deviations are calculated using 50 forward simulated solutions, and the number of iterations of VFGL is set to 200. The initial choice of parametric form has the best KKT deviation and solution error values. Table 7 shows that both the first stage solution error and objective value are inversely related to KKT deviation in the proposed parametric forms, which verifies the use of KKT deviation for the performance measure of parametric family choice.

*Table 7. Comparison of KKT deviation, first stage solution error, and objective value between various parametric forms for production planning problem (\* indicates the baseline case for the comparison)*

| Parametric form | KKT deviation | First stage solution error | Objective value |
|---|---|---|---|
| $*s_{t-1} \cdot \theta_t^1 + \sum_{i=1}^{3} e^{-\theta_{t,i}^2 s_{t-1,i}}$ | 2.53 (0.32) | 0.01 (0.01) | 178.40 (0.06) |
| $s_{t-1} \cdot \theta_t^1 + s_{t-1}^2 \cdot \theta_t^2$ | 49.14 (0.32) | 0.52 (0.00) | 180.19 (0.06) |
| $s_{t-1} \cdot \theta_t^1 + e^{s_{t-1} \cdot \theta_t^2}$ | 98.33 (0.32) | 1.36 (0.07) | 188.74 (0.06) |
| $s_{t-1} \cdot \theta_t^1$ | 110.30 (0.32) | 1.41 (0.000) | 188.74 (0.06) |
| 0 | 137.6 (0.32) | 1.41 (0.000) | 188.40 (0.06) |

The value of the parameter is slightly perturbed to check whether the same parametric form of the value function is still reusable. In particular, the maximum amount of resources $r_t$ is perturbed using the initial parametric form, and the result is summarized in Table 8. The result corresponding to the original experimental value ($r_t = 10$) is marked in the table. The parametric form of the experiment is applicable to the perturbation of the original problem. The parametric form of the experiment gives an accurate first stage solution and objective values for the upward perturbation of the maximum resource. However, downward perturbation of the maximum resource with a magnitude greater than 0.5 results in an erroneous first stage solution. Despite the error in the first stage solution, the objective value for the VFGL solution seems to remain accurate, which implies that either the solution from VFGL is optimal or the suboptimality of the VFGL solution is not substantial.

*Table 8. Comparison of perturbed problems with a fixed parametric value function form for production planning problem (\* indicates the baseline case for the comparison)*

| Maximum resource | VFGL | | | | | MSP | | | |
|---|---|---|---|---|---|---|---|---|---|
| | KKT deviation | Objective | Product 1 | Product 2 | Product 3 | Objective | Product 1 | Product 2 | Product 3 |
| 8 | 1.22 | 227.63 | 1.63 | 0.33 | 1.07 | 227.22 | 2.00 | 0.00 | 1.20 |
| 8.5 | 1.85 | 215.83 | 1.61 | 0.09 | 1.07 | 214.65 | 2.50 | 0.00 | 1.10 |
| 9 | 0.98 | 202.38 | 1.58 | 0.00 | 1.02 | 202.44 | 2.00 | 0.00 | 1.00 |
| 9.5 | 1.35 | 190.08 | 1.49 | 0.00 | 1.01 | 190.39 | 1.50 | 0.00 | 1.00 |
| *10 | 2.50 | 177.51 | 0.99 | 0.00 | 1.01 | 178.33 | 1.00 | 0.00 | 1.00 |
| 10.5 | 2.75 | 167.59 | 0.49 | 0.00 | 1.01 | 167.50 | 0.50 | 0.00 | 1.00 |
| 11 | 2.55 | 156.63 | 0.00 | 0.00 | 1.00 | 156.67 | 0.00 | 0.00 | 1.00 |

| | | | | | | | | |
|---|---|---|---|---|---|---|---|---|
| 11.5 | 2.24 | 147.74 | 0.00 | 0.00 | 0.90 | 147.67 | 0.00 | 0.00 | 0.90 |
| 12 | 1.78 | 138.50 | 0.00 | 0.00 | 0.81 | 138.67 | 0.00 | 0.00 | 0.80 |
| 12.5 | 2.00 | 130.38 | 0.00 | 0.00 | 0.71 | 129.67 | 0.00 | 0.00 | 0.70 |
| 13 | 2.03 | 121.94 | 0.00 | 0.00 | 0.61 | 120.67 | 0.00 | 0.00 | 0.60 |
| 13.5 | 2.17 | 111.65 | 0.00 | 0.00 | 0.51 | 111.67 | 0.00 | 0.00 | 0.50 |
| 14 | 4.16 | 102.90 | 0.00 | 0.00 | 0.40 | 102.67 | 0.00 | 0.00 | 0.40 |
| 14.5 | 2.76 | 94.51 | 0.00 | 0.00 | 0.20 | 94.00 | 0.00 | 0.00 | 0.20 |
| 15 | 4.70 | 85.08 | 0.00 | 0.00 | 0.00 | 85.33 | 0.00 | 0.00 | 0.00 |

## 4.3 Hydrothermal generation

In this example, the aim is to find an optimal electricity generation level for hydro and thermal plants while maintaining a desirable reservoir level. This example is a slightly modified and simplified version of the example given by Guigues (2014). For each stage $t$, the deterministic demand of electricity has to be met. The electricity system contains hydro plants and thermal plants for generation. Hydro plants can produce electricity at a lower price, but the resource is limited by the reservoir level. In a hydro plant, the water inflow to the reservoir is uncertain. It is assumed that thermal plants have unlimited capacity in electricity generation but require a higher generation cost. In addition to the cost of electricity generation, there is a disutility cost with respect to each stage of the final water reservoir level. The disutility reflects the environmental concern while draining the reservoir. The objective is to minimize the sum of the expected generation cost and disutility. For the case study, a model with seven stages is considered. We assume that the water inflow to the reservoir follows independent and identical normal distributions over different stages. The decision variables, parameters with their values for this example are described in the Table 9.

*Table 9. Decision variables and parameters for hydrothermal generation problem*

| Decision variables | Description | |
|---|---|---|
| $r_t^{init}$ | Water reservoir level in the beginning of stage $t$ | |
| $r_t^{final}$ | Water reservoir level in the end of stage $t$ | |
| $W_t$ | Hydro electricity generation level at stage $t$ | |
| $H_t$ | Thermal electricity generation level at stage $t$ | |
| **Parameters** | **Description** | **Value** |
| $r_0^{init}$ | Initial water reservoir level | 40 |
| $c_t^W$ | Cost of hydro electricity production per unit at stage $t$ | 2 |
| $c_t^H$ | Cost of thermal electricity production per unit at stage $t$ | 7 |
| $d_t$ | Electricity demand at stage $t$ | 20 |
| $a_t$ | Reservoir level utility coefficient | 0.1 |
| $b_t$ | Reservoir level utility scaling constant | 5 |
| $I_t$ | Water inflow to reservoir in the beginning of stage $t$ | Normally distributed with mean 20, standard deviation 5 |

The total cost is as follows.

$$\mathbb{E}\left[\sum_{t=1}^{T} c_t^W W_t + c_t^H H_t + e^{-a_t r_t^{final}}\right]$$

The stagewise subproblems are given below. Here, the third term of the objective function is a disutility on the reservoir level, which can be considered a negative of scaled exponential utility.

### *Stage 1 subproblem*

minimize $\quad c_1^W W_1 + c_1^H H_1 + e^{-a_1 r_1^{final} + b_1} + V_2(r_1^{final})$

subject to
| | | |
|---|---|---|
| $r_1^{init} = r_0$ | | Initial reservoir |
| $r_1^{final} = r_1^{init} - W_1$ | | Reservoir balance |
| $W_1 + H_1 \geq d_1$ | | Demand |
| $r_1^{final}, W_1, H_1 \geq 0$ | | Non-negativity |

### *Stage t subproblem* $(t = 2, ..., T-1)$

minimize $\quad c_t^W W_t + c_t^H H_t + e^{-a_t r_t^{final} + b_t} + V_{t+1}(r_t^{final})$

subject to
| | | |
|---|---|---|
| $r_t^{init} = r_{t-1}^{final} + I_t$ | | Initial reservoir |
| $r_t^{final} = r_t^{init} - W_t$ | | Reservoir balance |
| $W_t + H_t \geq d_t$ | | Demand |
| $r_t^{final}, W_t, H_t \geq 0$ | | Non-negativity |

### *Stage T subproblem*

minimize $\quad c_T^W W_T + c_T^H H_T + e^{-a_T r_T^{final} + b_T}$

subject to
| | | |
|---|---|---|
| $r_T^{init} = r_{T-1}^{final} + I_t$ | | Initial reservoir |
| $r_T^{final} = r_T^{init} - W_T$ | | Reservoir balance |
| $W_T + H_T \geq d_T$ | | Demand |
| $r_T^{final}, W_T, H_T \geq 0$ | | Non-negativity |

Next, the parametric value function form $\hat{V}_t$ must be determined. The parametric utility function is set as follows.

$$\hat{V}_t\left(r_{t-1}^{final}; \theta_t = (\theta_t^1, \theta_t^2)\right) = \theta_t^1 r_{t-1}^{final} + e^{-\theta_t^2 r_{t-1}^{final} + b_t}, t = 2, ..., T$$

where the first linear term is added to capture the cost reduction by increasing the portion of hydro plant generation, and the second term is added to capture the decreased disutility from the increased reservoir level. The second term is the indefinite integral form of the stagewise objective function proposed in Remark 1. The initial parameter values are set as follows.

$$(\theta_t^1, \theta_t^2) = (-1, -1), \quad t = 2, ..., T$$

The scaling term of the parametric value function is fixed by $b_t$, which is the scaling constant of the utility function from the original problem. Therefore, it is advisable to try first a parametric value function with minimal parameter flexibility and then increase the parameter flexibility if the result is not satisfactory.

*Table 10. Comparison of objective value, first stage solution, and computation time between algorithms for hydrothermal generation problem*

| Algorithm | Objective value | Hydro plant generation | Thermal plant generation | Computation time (s) |
|---|---|---|---|---|
| **Optimal (MSP)** | 397 (0.1) | 9.90 (0.05) | 10.10 (0.05) | - |
| **SDDP** | 397 (2) | 9.78 (0.41) | 10.22 (0.41) | 978 (0.20) |
| **VFGL** | 400 (1) | 10.45 (0.15) | 9.55 (0.15) | 337 (3.15) |

The experimental results of MSP, SDDP, and VFGL are summarized in Table 10. It presents the first stage solution of the true MSP model, and the averaged first stage solutions from 20 generated instances of SDDP and VFGL. The table shows that the MSP and SDDP give similar solutions, whereas the result from VFGL is slightly different. A slight suboptimality of VFGL in the objective value compared to those of MSP and SDDP is also observed.

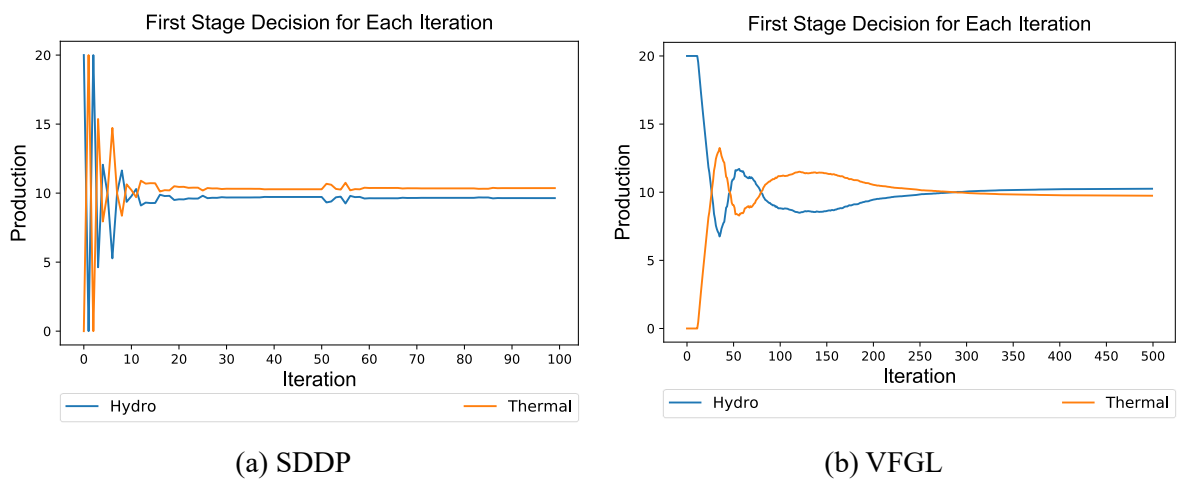

(a) SDDP      (b) VFGL

**Figure 5**. Comparison of evolution of the first stage solution between SDDP and VFGL for hydrothermal generation problem

Figure 5 displays a sample of the first-stage solution convergence for SDDP and VFGL among 20 trials. It shows that SDDP and VFGL seem to converge after 100 and 500 iterations, respectively. Table 10 and Figure 5 indicate that the VFGL tends to underappreciate slightly the thermal electric generation compared with SDDP and MSP. The applicability of VFGL may depend on the tolerance level.

However, two interesting properties of VFGL can be found in this experiment. First, the evolution of solutions from VFGL tends to fluctuate less than that of SDDP. It is true that the solutions from the SDDP algorithm quickly converged after approximately 10 iterations. However, during the early learning phase, the solutions fluctuated drastically with every iteration. This implies that it is extremely dangerous to use the result from the SDDP algorithm that has not sufficiently converged because it might yield extremely different results after more iterations. However, the degree of solution fluctuation is much smaller for the VFGL algorithm. There may be a gradual update of the value function for the stability of the solution convergence of VFGL. In other words, a small update of the value function parameters is smoother and less drastic than adding a linear cut to the value function approximation. Second, the computation time for VFGL is approximately one third that of SDDP, even with a five times greater number of iterations. Similar to the first two examples, the computation time of SDDP increases in every iteration, while that of VFGL remains relatively consistent (see Appendix G). Therefore, the VFGL algorithm would be more computationally robust for large-scale problems that require a very large number of iterations.

This example is solved with different parametric forms of the value function, and the results are reported in Table 11. A lower KKT deviation leads to less solution error. However, the first four parametric forms give a comparable objective. Therefore, the suboptimality of solutions from the first four parametric forms is likely to be small. The quadratic parametric function includes a linear term, implying that its model capacity is greater than that of the linear function. However, there are a higher KKT deviation and solution error for the quadratic form compared with the linear one. This implies that the trained parameters for the quadratic function fell into local optimality. The suboptimality of the parameters occurs owing to the nonconvexity of the loss function with respect to the parameters. Nevertheless, the quadratic form gives a better objective value. For this reason, it is suspected that there are a flat objective value and approximation error.

*Table 11. Comparison of KKT deviation, first stage solution error, and objective value between various parametric forms for hydrothermal generation problem (\* indicates the baseline case for the comparison)*

| Parametric form | KKT deviation | First stage solution error | Objective value |
| --- | --- | --- | --- |
| *$\theta_t^1 r_{t-1}^{final} + e^{-\theta_t^2 r_{t-1}^{final} + b_t}$ | 3.56 (1.27) | 0.8 (0.8) | 400 (0.8) |
| $\theta_t^1 r_{t-1}^{final}$ | 3.65 (1.27) | 1.1 (0.7) | 400.8 (0.1) |
| $\theta_t^1 r_{t-1}^{final} + \theta_t^2 (r_{t-1}^{final})^2$ | 3.50 (1.68) | 1.3 (0.7) | 399.0 (0.1) |
| $\theta_t^1 r_{t-1}^{final} - \log(\theta_t^2 r_{t-1}^{final})$ | 5.98 (1.06) | 4.2 (1.7) | 400.1 (0.4) |

| | | 14.02 (2.07) | 14.2 (0.0) | 451.1 (0.0) |
|---|---|---|---|---|
| 0 | | | | |

The original problem parameters are perturbed to determine whether the parametric value function is still applicable. In particular, the costs of hydro and thermal generation vary. For most of the perturbations in Table 12, VFGL gives a slightly suboptimal objective value with a noticeable difference in the first stage solution compared with the results of MSP. The suboptimality of VFGL seems to be greater for the thermal generation cost perturbation than for the hydro generation cost perturbation. Depending on the tolerance level, the suboptimality of VFGL with a given parametric form might be acceptable or unacceptable. To make the VFGL solution more accurate for the perturbed problems, more suitable parametric families of the value function should be explored.

*Table 12. Comparison of perturbed problems with a fixed parametric value function form for hydrothermal generation problem (\* indicates the baseline case for the comparison)*

| Parameter | | VFGL | | | | MSP | | |
|---|---|---|---|---|---|---|---|---|
| $c_t^W$ | $c_t^H$ | KKT deviation | Objective value | Hydro plant generation | Thermal plant generation | Objective value | Hydro plant generation | Thermal plant generation |
| 2 | 5 | 2.24 | 365.00 | 5.06 | 14.94 | 370.74 | 5.39 | 14.61 |
| 2 | 5.5 | 2.68 | 374.01 | 6.26 | 13.74 | 378.38 | 7.53 | 12.47 |
| 2 | 6 | 2.45 | 382.01 | 7.68 | 12.32 | 378.99 | 8.36 | 11.64 |
| 2 | 6.5 | 3.24 | 389.62 | 8.91 | 11.09 | 392.04 | 9.8 | 10.2 |
| *2 | *7 | 3.66 | 396.64 | 9.90 | 10.10 | 400.39 | 10.45 | 9.55 |
| 2 | 7.5 | 3.6 | 402.68 | 10.91 | 9.09 | 404.15 | 10.88 | 9.12 |
| 2 | 8 | 4.52 | 408.99 | 11.61 | 8.39 | 417.47 | 13.18 | 6.82 |
| 3 | 7 | 2.66 | 521.85 | 7.67 | 12.33 | 523.37 | 8.68 | 11.32 |
| 3.5 | 7 | 2.44 | 583.90 | 6.49 | 13.51 | 585.67 | 6.89 | 13.11 |
| 4 | 7 | 1.87 | 645.08 | 4.88 | 15.12 | 644.56 | 5.11 | 14.89 |
| 4.5 | 7 | 1.89 | 705.74 | 2.78 | 17.22 | 711.22 | 3.91 | 16.09 |
| 5 | 7 | 1.35 | 764.90 | 0.81 | 19.19 | 765.76 | 1.54 | 18.46 |

## 4.4 Discussion

Throughout the experiment, the VFGL algorithm successfully solved the problems with an accuracy comparable to that of the SDDP, given the right parametric form for the value function approximation. A key characteristic of the VFGL is an almost constant computation time per iteration, which is shown to be a clear computational advantage compared to SDDP.

As a critical drawback of the algorithm, the parametric form of the value function must be determined prior to the problem solving. However, this is not a specific problem with our algorithm. Ghadimi, Perkins, and Powell (2020) noted that the choice of appropriate parameterization is often considered an art of modeling for engineering problems, particularly within the machine learning domain. Indeed, it was shown to be a vital part of the proposed VFGL algorithm, as such choices had a significant impact on the solution quality. Therefore, we introduced a measure for evaluating the fitness of a given parametric form called the KKT deviation, which was proven to be useful for evaluating different parameterizations. One rule of thumb for the choice of parametric form would be to start from the same functional form of the objective function or the indefinite integral form of the objective function.

## 5. Conclusion

A novel learning-based stagewise decomposition algorithm, VFGL, was proposed to solve large-scale multistage stochastic programs approximately. The algorithm approximates the value functions by a convex parametric form of functions, which is distinguished from the Benders decomposition based algorithms that use piecewise linear approximation. Three illustrative examples demonstrated the computational advantages of VFGL compared with the deterministic equivalent formulation approach and widely used stagewise decomposition algorithm called SDDP. The choice of appropriate parametric form would require some effort; thus, the KKT deviation was proposed to measure the suitability of different parametric forms. Finally, numerical experiments indicated that VFGL can recycle the same parametric form for the same optimization problem with slightly different parameters, which can be useful for real-world applications (Bertsimas & Stellato, 2020).

Future studies on VFGL should focus on finding an appropriate parametric form of the value function. One possible approach would be to construct a general parametric function approximation with a large modeling capacity, similar to artificial neural networks. Another possible method is directly parametrizing the gradient of the value function, which may greatly increase the modeling capacity.

**Acknowledgement**

This research was supported by Basic Science Research Program through the National Research Foundation of Korea (NRF) funded by the Ministry of Science and ICT (NRF-2020R1A2C101067711 and NRF-2019R1C1C1010456).